\documentstyle[11pt,figure]{article}

\begin{document}

\font\bbbld=msbm10 scaled\magstep1
\newcommand{\bfR}{\hbox{\bbbld R}}
\newcommand{\bfC}{\hbox{\bbbld C}}
\newcommand{\bfZ}{\hbox{\bbbld Z}}
\newcommand{\bfH}{\hbox{\bbbld H}}
\newcommand{\bfQ}{\hbox{\bbbld Q}}
\newcommand{\bfN}{\hbox{\bbbld N}}
\newcommand{\bfP}{\hbox{\bbbld P}}
\newcommand{\bfT}{\hbox{\bbbld T}}
\def\Sym{\mathop{\rm Sym}}
\newcommand{\halo}[1]{\Int(#1)}
\def\Int{\mathop{\rm Int}}
\def\Re{\mathop{\rm Re}}
\def\Im{\mathop{\rm Im}}
\newcommand{\union}{\cup}
\newcommand{\goesto}{\rightarrow}
\newcommand{\bdy}{\partial}
\newcommand{\n}{\noindent}
\newcommand{\p}{\hspace*{\parindent}}

\newtheorem{theorem}{Theorem}[section]
\newtheorem{assertion}{Assertion}[section]
\newtheorem{proposition}{Proposition}[section]
\newtheorem{lemma}{Lemma}[section]
\newtheorem{definition}{Definition}[section]
\newtheorem{claim}{Claim}[section]
\newtheorem{corollary}{Corollary}[section]
\newtheorem{observation}{Observation}[section]
\newtheorem{conjecture}{Conjecture}[section]
\newtheorem{question}{Question}[section]
\newtheorem{example}{Example}[section]

\newbox\qedbox
\setbox\qedbox=\hbox{$\Box$}
\newenvironment{proof}{\smallskip\noindent{\bf Proof.}\hskip \labelsep}%
                        {\hfill\penalty10000\copy\qedbox\par\medskip}
\newenvironment{remark}{\smallskip\noindent{\bf Remark.}\hskip \labelsep}%
                        {\hfill\penalty10000\copy\qedbox\par\medskip}
\newenvironment{remark2}{\smallskip\noindent{\bf Remark 2.}\hskip \labelsep}%
                        {\hfill\penalty10000\copy\qedbox\par\medskip}
\newenvironment{proofspec}[1]%
                      {\smallskip\noindent{\bf Proof of Theorem 1.1.}
                        \hskip \labelsep}%
                        {\nobreak\hfill\hfill\nobreak\copy\qedbox\par\medskip}
\newenvironment{proofspec2}[1]%
                      {\smallskip\noindent{\bf Proof of Theorem 1.2.}
                        \hskip \labelsep}%
                        {\nobreak\hfill\hfill\nobreak\copy\qedbox\par\medskip}
\newenvironment{acknowledgements}{\smallskip\noindent{\bf Acknowledgements.}%
        \hskip\labelsep}{}

\setlength{\baselineskip}{1.4\baselineskip}

\title{On Embeddedness of Area-Minimizing Disks, 
and an Application to Constructing Complete Minimal Surfaces}
\author{Wayne Rossman}
\vspace{-0.2in} 
\maketitle

\begin{abstract}
Let $\alpha$ be a polygonal Jordan curve in $\bfR^3$.  
We show that if $\alpha$ satisfies certain 
conditions, then the least-area Douglas-Rad\'{o} disk 
in $\bfR^3$ with boundary $\alpha$ is unique and is a 
smooth graph.  
As our conditions on $\alpha$ are not included amongst 
previously known conditions for
embeddedness, we are enlarging the set of
Jordan curves in $\bfR^3$ which are 
known to be spanned by an embedded least-area disk.

As an application, 
we consider the conjugate surface construction method for minimal surfaces.
With our result we can apply this method 
to a wider range of complete catenoid-ended minimal
surfaces in $\bfR^3$.
\footnote{This research was
supported by a fellowship from the Japan
Society for the Promotion of Science.}
\end{abstract}

\section{Introduction}

Much investigation has been made on the Plateau problem, i.e. 
to show that any rectifiable Jordan curve in $\bfR^3$ 
bounds a minimal surface of least area.  The first results were 
by Douglas and Rad\'{o} in the early 1930's, 
when they proved existence of a smooth least-area disk 
for any given boundary
curve \cite{Os}.  This disk is often called the Douglas-Rad\'{o}
solution.  
Osserman \cite{Os} later showed that the Douglas-Rad\'{o} solution has no
true branch points in its interior, and Gulliver \cite{Gu}
showed that it also has no false branch points in its interior.  
Hildebrandt showed regularity at the boundary of the 
Douglas-Rad\'{o} solution wherever its boundary is real-analytic.  
Then, along real-analytic portions of the boundary, 
Gulliver and Lesley \cite{GuLe} showed nonexistence of branch points.  
Putting all of these results together, we have the following theorem.  

\begin{theorem}
Let $\alpha$ be a rectifiable Jordan curve in $\bfR^3$.  Then there
exists a map
\[ h: D \rightarrow \bfR^3 \; ,\]
where $D$ is the closed unit disk in $\bfR^2$, satisfying
\begin{enumerate}
\item $h$ is continuous in $D$;
\item $h$ maps the boundary of $D$, $\partial D$, bijectively to $\alpha$;
\item $h \in C^\infty$ (in fact, $h$ is harmonic) in the interior of
  $D$, and is a regular conformal minimal immersion in the 
  interior of $D$; 
\item the image of $D$ under $h$ has minimum area among all maps $D
\rightarrow \bfR^3$ which are piecewise smooth in the interior and
satisfy the conditions 1 and 2 above.  
\item If $\sigma$ is a closed subarc of $\partial D$ that is mapped
  by $h$ into the interior of some real-analytic subarc $\gamma$ of
  $\alpha$, then $h$ can be analytically
  continued across $\sigma$ (as a minimal surface), and $h$ has no 
  branch points on $\sigma$.  
\end{enumerate}
\end{theorem}

In the case that the $\gamma$ above 
is a straight line or a planar geodesic, an
even stronger conclusion is known, 
and is called the Schwarz reflection principle (\cite{Ka2}, Sec. 1.3.2):

\begin{theorem}
The union of a minimal surface with its reflection (resp.
rotation by 180 degrees) 
across (resp. about) a plane containing a boundary planar
geodesic (resp. a line segment in the boundary of the 
surface) is a smooth minimal surface.
\end{theorem}

The results above go far to solve the problem of 
existence and regularity of 
least-area disks for a given curve, as well as to show 
nonexistence of branch points in least-area disks.  
However, the question of embeddedness is only partly 
answered.  Almgren and Thurston showed that there exist 
unknotted Jordan curves that cannot bound any embedded minimal disk
\cite{MeYa1}.  It seems 
difficult to find conditions on a curve that imply 
its Douglas-Rad\'{o} solutions are 
embedded, but some partial results have been found.  
Rad\'o proved in 1932 \cite{MeYa2} that if 
$\alpha$ is an embedded rectifiable curve in $\bfR^3$ whose vertical 
projection to the $x_1x_2$-plane (or a central projection 
from a certain point) is one-to-one and convex, 
then the Douglas-Rad\'{o}
solution is the unique least-area surface bounded by $\alpha$ and 
is a graph over the $x_1x_2$-plane (or a graph with respect to 
central projection).  
Meeks and Yau \cite{MeYa1} generalized this to the case that 
$\alpha$ is extremal, i.e. lies on the boundary of its convex 
hull.  They showed that 
any Douglas-Rad\'{o} solution for an extremal curve $\alpha$ is embedded.
They later generalized this to show the same conclusion even when 
$\alpha$ only lies in the
boundary $\partial \hat{M}$ 
of a closed region $\hat{M} \subseteq \bfR^3$ such that 
$\partial \hat{M}$ has nonnegative mean curvature with respect to
the interior of $\hat{M}$ \cite{MeYa2}.

We shall show (Theorem 2.1) that for certain types of
polygonal Jordan curves in $\bfR^3$, the Douglas-Rad\'{o} solution is
an embedded graph.  We can apply Theorem 2.1 in some cases where the
results of Meeks and Yau do not apply.  
The original motivation for considering these
types of polygonal Jordan curves is their usefulness in 
the conjugate surface construction
method for minimal surfaces in $\bfR^3$ 
(\cite{BeRo}, \cite{Ka1}, \cite{Ka2}, \cite{Ka3}, \cite{Ka4}, 
\cite{Ro}).  Some examples of this construction are shown in 
Section 5.  

Theorem 2.1 allows us to extend the
conjugate surface construction to more cases.  
The strategy is roughly as follows:  We wish to prove existence of 
complete 
catenoid-ended minimal surfaces $M$ with symmetry, where the symmetries
of $M$ are generated by a discrete set of reflections in 
$\bfR^3$.  We consider the smallest portion of $M$ that will generate the
entire surface under the action of the symmetry group, and we call
it the {\em fundamental piece} of $M$.  We 
choose the fundamental piece so that it is bounded by planar 
geodesics.  It is then enough to show existence of the fundamental
piece only, since the entire surface $M$ can be produced from the fundamental
piece by reflection (Theorem 1.2).  Furthermore, we can 
show existence of the fundamental piece by showing existence of
the conjugate surface $M^\prime$ to the fundamental piece.  
(We define the conjugate surface in Section 3.)  
The advantage of considering the conjugate
surface $M^\prime$ is that it is bounded by straight lines.  
We prove the existence of $M^\prime$
by showing it exists as the limit of a sequence of compact
embedded stable minimal disks $M_i$ bounded by Jordan 
polygonal curves $\alpha_i$.  

We will use the term stable in the following sense:
Minimal surfaces are critical for the first variation formula.  
A minimal surface $\cal S$ (possibly with boundary $\partial 
{\cal S}$) is {\em stable} if the
second derivative of area is nonnegative at $\cal S$ for all
smooth variations of the surface with compact support (and fixing 
$\partial {\cal S}$).  

So the first step is to demonstrate the existence of $M_i$ 
bounded by $\alpha_i$.  For the minimal surfaces $M$ we are
considering, $\alpha_i$ can be chosen to satisfy all the conditions of
Theorem 2.1.  
Thus the Douglas-Rad\'{o} surfaces $M_i$ for $\alpha_i$ 
are smooth graphs in $\bfR^3$.  In particular, 
$M_i$ are smoothly embedded and stable.  Once we have stability, we
can show that $\{M_i\}_{i=1}^\infty$ has a convergent subsequence
(Lemma 4.1).  $M^\prime$ is the limit surface.  
(The question of minimal graphs over unbounded planar domains has been
investigated in \cite{EaRo}, \cite{BeRo}, and \cite{Ro}.)  
We then show that $M^\prime$ is 
connected in the cases we consider.

In the case that $M$ may have some unwanted periodicity,
we need to 
show that $M^\prime$ can be constructed so that $M$ does not have this
periodicity.  Lemma 4.2 is useful for this.  

The author thanks Miyuki Koiso, 
Shin Nayatani, Leon Simon, the referee, and members of
G.A.N.G. for helpful suggestions.

\section{The Main Result}

\begin{theorem}
Let $\alpha = \ell_1 \cup ... \cup \ell_m$ be a 
closed embedded polygonal curve in $\bfR^3$ consisting of 
straight line segments $\ell_i$ and vertices 
$\ell_1 \cap \ell_2, ..., \ell_{m-1} \cap \ell_m, \ell_m \cap \ell_1$.  
Let $P$ be a polygonal region in the $x_1x_2$-plane $\{x_3 = 0 \}$ bounded 
by the polygon $\partial P = \rho_1 \cup ... \cup \rho_{m-1}$ consisting 
of edges $\rho_i$ and vertices $\rho_1 \cap \rho_2,... , 
\rho_{m-2} \cap \rho_{m-1}, \rho_{m-1} \cap \rho_1$.  
Suppose the following:
\begin{enumerate}
\item $\ell_m$ is vertical; and for each $i = 1,...,m-1$, $\ell_i$ is not 
  vertical.  
\item $\ell_{m-1}$ and $\ell_1$ are horizontal, and $\alpha$ lies entirely 
  between the two horizontal planes containing $\ell_{m-1}$ and $\ell_1$.  
  That is, there exist $a,b \in \bfR$ such that 
  $\ell_{m-1} \subset \{x_3 = a \}$ and 
  $\ell_{1} \subset \{x_3 = b \}$ and $\alpha \subset \{
  \mbox{min}(a,b) \leq x_3 \leq \mbox{max}(a,b) \}$.  
\item Denoting the boundary of the convex hull
  of $P$ in the $x_1x_2$-plane by $\partial \mbox{Conv}(P)$, 
  we have $\partial \mbox{Conv}(P) \cap \partial P = \rho_2 \cup ... 
  \cup \rho_{m-2}$.  
\item Each $\ell_i$, $i=1,...,m-1$ is mapped bijectively to $\rho_i$ by 
  the vertical projection 
  ${\cal P}: (x_1,x_2,x_3) \to (x_1,x_2,0)$, and ${\cal P}(\ell_m) = 
  \rho_{m-1} \cap \rho_1$.  
\end{enumerate}
Then the 
Douglas-Rad\'{o} solution with boundary $\alpha$ 
is unique and embedded, and its interior is a graph over the interior
of $P$. 
\end{theorem}

\begin{remark}
It is clear from the proof below that this theorem could be 
generalized somewhat.  For example, we could easily adapt the proof to 
include cases where $\ell_{m-1}$ and $\ell_1$ are not horizontal, or where 
$\alpha$ has portions that are not 
polygonal.  However, as the statement above is sufficient for the 
applications in Section 5, for simplicity we do not consider any 
generalizations here.
\end{remark}

\begin{figure}
\hspace{0.5in}
\unitlength=1.0pt
\begin{picture}(356.00,175.00)(113.00,653.00)
\put(280.00,791.00){\line(0,-1){61.00}}
\put(280.00,730.00){\line(5,-1){48.00}}
\put(328.00,721.00){\line(-1,0){192.00}}
\put(136.00,721.00){\line(3,5){64.00}}
\put(200.00,828.00){\line(1,0){124.00}}
\put(324.00,828.00){\line(-6,-5){44.00}}
\put(280.00,669.00){\line(5,-1){48.00}}
\put(328.00,660.00){\line(-1,0){192.00}}
\put(136.00,660.00){\line(5,3){64.00}}
\put(200.00,699.00){\line(1,0){125.00}}
\put(325.00,699.00){\line(-3,-2){45.00}}
\put(420.00,740.00){\vector(0,1){40.00}}
\put(420.00,740.00){\vector(1,0){40.00}}
\put(420.00,740.00){\vector(-1,-4){7.00}}
\put(354.00,658.00){\makebox(0,0)[cc]{\small $(\lambda_2,\beta_2,0)$}}
\put(262.00,671.00){\makebox(0,0)[cc]{\small $(0,0,0)$}}
\put(354.00,698.00){\makebox(0,0)[cc]{\small $(-\lambda_1,\beta_1,0)$}}
\put(203.00,706.00){\makebox(0,0)[cc]{\small $(-\lambda_1,-\gamma,0)$}}
\put(136.00,653.00){\makebox(0,0)[cc]{\small $(\lambda_2,-\gamma,0)$}}
\put(352.00,723.00){\makebox(0,0)[cc]{\small $(\lambda_2,\beta_2,1)$}}
\put(262.00,732.00){\makebox(0,0)[cc]{\small $(0,0,1)$}}
\put(262.00,792.00){\makebox(0,0)[cc]{\small $(0,0,2)$}}
\put(354.00,826.00){\makebox(0,0)[cc]{\small $(-\lambda_1,\beta_1,2)$}}
\put(168.00,828.00){\makebox(0,0)[cc]{\small $(-\lambda_1,-\gamma,2)$}}
\put(133.00,713.00){\makebox(0,0)[cc]{\small $(\lambda_2,-\gamma,1)$}}
\put(307.00,807.00){\makebox(0,0)[cc]{\tiny $\ell_1$}}
\put(250.00,824.00){\makebox(0,0)[cc]{\tiny $\ell_2$}}
\put(155.00,763.00){\makebox(0,0)[cc]{\tiny $\ell_3$}}
\put(210.00,725.00){\makebox(0,0)[cc]{\tiny $\ell_4$}}
\put(305.00,729.00){\makebox(0,0)[cc]{\tiny $\ell_5$}}
\put(276.00,765.00){\makebox(0,0)[cc]{\tiny $\ell_6$}}
\put(305.00,680.00){\makebox(0,0)[cc]{\tiny $\rho_1$}}
\put(272.00,702.00){\makebox(0,0)[cc]{\tiny $\rho_2$}}
\put(164.00,681.00){\makebox(0,0)[cc]{\tiny $\rho_3$}}
\put(220.00,663.00){\makebox(0,0)[cc]{\tiny $\rho_4$}}
\put(302.00,668.00){\makebox(0,0)[cc]{\tiny $\rho_5$}}
\put(225.00,680.00){\makebox(0,0)[cc]{$P$}}
\put(180.00,777.00){\makebox(0,0)[cc]{$\alpha$}}
\put(413.00,704.00){\makebox(0,0)[cc]{$x_1$}}
\put(469.00,740.00){\makebox(0,0)[cc]{$x_2$}}
\put(419.00,788.00){\makebox(0,0)[cc]{$x_3$}}
\end{picture}
        \caption{A curve $\alpha$ satisfying all the conditions of
          Theorem 2.1.}
\end{figure}
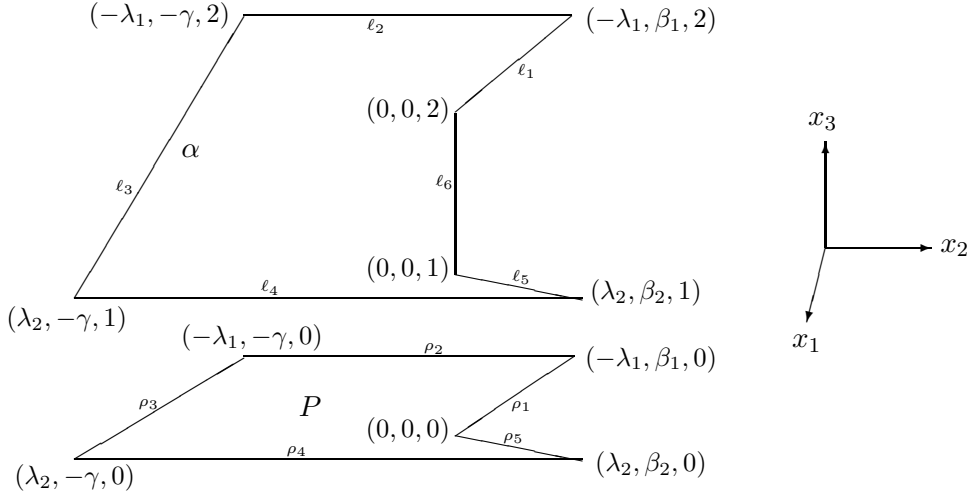

\begin{example}
{\em Let $\lambda_1,\lambda_2,\beta_1,\beta_2,\gamma$ be 
any positive numbers.  
Let $\alpha$ be the polygonal curve from $(0,0,1)$ to $(0,0,2)$ to
$(-\lambda_1,\beta_1,2)$ to $(-\lambda_1,-\gamma,2)$ to
$(\lambda_2,-\gamma,1)$ to $(\lambda_2,\beta_2,1)$ and back to
$(0,0,1)$.  Let $\partial P$ be the 5-gon in the $x_1x_2$-plane 
with vertices $(0,0,0)$, 
$(\lambda_2,\beta_2,0)$, $(-\lambda_1,\beta_1,0)$, 
$(-\lambda_1,-\gamma,0)$, and $(\lambda_2,-\gamma,0)$, so that 
$\partial P = {\cal P}(\alpha)$ (see Figure 1).  
By Theorem 2.1, the Douglas-Rad\'{o} solution for 
$\alpha$ is unique and is a graph over $P$.  

In the case $\lambda_1 = \lambda_2$ and $\beta_1 = \beta_2$,
it was already known that $\alpha$ bounds a smoothly embedded minimal disk
that is stable in $\bfR^3$.  Consider the polygonal curve
$\tilde{\alpha}$ from $(0,0,1)$ to $(\lambda_2,\beta_2,1)$ to
$(\lambda_2,-\gamma,1)$ to $(0,-\gamma,3/2)$ to
$(0,0,3/2)$ and back to $(0,0,1)$.  The least-area surface $\tilde{M}$
spanning $\tilde{\alpha}$ is unique and is a graph over the region 
$P \cap \{x_1 \geq 0\}$, by 
Nitsche's theorem \cite{BeRo}, \cite{Ro}.  Let Rot$: \bfR^3 
\rightarrow \bfR^3$ be 
rotation by 180 degrees about the line
through $(0,-\gamma,3/2)$ and $(0,0,3/2)$.  
Then $\tilde{M} \cup $Rot($\tilde{M}$)
is a smoothly embedded minimal graph with boundary $\alpha$, by
Theorem 1.2.  Since the image of the Gauss map on 
$\tilde{M} \cup $Rot($\tilde{M}$) is contained in
a hemisphere, $\tilde{M} \cup $Rot($\tilde{M}$) is stable \cite{BdC}.

However, Theorem 2.1 shows that $\tilde{M} \cup $Rot($\tilde{M}$) is
also the unique least-area surface with boundary $\alpha$.  In
fact, Example 2.1 shows existence of a unique least-area surface of disk 
type with
boundary $\alpha$ in the nonsymmetric cases $\lambda_1 \neq \lambda_2$ or 
$\beta_1 \neq \beta_2$ as well, where it was not previously
known if there were even stable embedded minimal disks with boundary
$\alpha$.  
}
\end{example}

\begin{remark}
Theorem 2.1 is not true without the fourth condition.  
For example, let 
$\alpha$ be the polygonal curve consisting of line segments from
$(0,0,0)$ to $(2,0,0)$ to $(2,\delta,0)$ to $(2,0,\delta)$ to 
$(0,0,\delta)$ to $(0,1,\delta)$ to 
$(\frac{1}{2},\epsilon,\delta)$ to $(\frac{1}{2},\epsilon,0)$ to 
$(0,1,0)$ 
and back to $(0,0,0)$.  Let $\partial 
P$ be the polygonal 5-gon in the $x_1x_2$-plane
with vertices $(0,0,0)$, $(2,0,0)$, $(2,\delta,0)$, $(0,1,0)$, and
$(\frac{1}{2},\epsilon,0)$, so that ${\cal P}(\alpha) \subset 
\partial P$.  
If $0 < \epsilon << \delta << 1$, then the least-area surface
bounded by $\alpha$ is not embedded, since its
interior will intersect $\alpha$ along the line segment from 
$(\frac{1}{2},\epsilon,0)$ to $(\frac{1}{2},\epsilon,\delta)$.  
\end{remark}

We now state two lemmas following from Theorem 4 and Lemmas 2 and 3 
of \cite{MeYa1}.  We use these two lemmas in the proof of Theorem 2.1.  
Let $B_\epsilon (p) := \{q \in \bfR^3 \; | \; \mbox{dist}(p,q) < \epsilon \}$.  

\begin{lemma}
If the self-intersection set 
$S(h) = \{p \in D \; |
\; \exists q \neq p \in D \mbox{ with } h(p) = h(q)\}$ is disjoint from 
$\partial D$, then $h$ is an embedding.  
\end{lemma}

\begin{lemma}
Let $g: D \rightarrow \bfR^3$ and $f: D \rightarrow \bfR^3$ be regular
minimal
embeddings that intersect at a point $p \in \bfR^3$ such that $p
\not\in g(\partial D) \cup f(\partial D)$.  Assume that the images of
$g$ and $f$ do not coincide in a neighborhood of $p$.  
Then for some small $\epsilon$,
the intersection set $f(D) \cap g(D) \cap B_\epsilon(p)$ 
consists of a finite number of curves through $p$ and the
intersection is transverse at points other than $p$.
The intersection set cannot be a point, and cannot contain a curve
with an endpoint in Int($B_\epsilon(p)$), and cannot have nonempty
interior.

In particular, this holds for the intersection of
a nonflat minimal immersion with any of its tangent planes.  
\end{lemma}

The proof of Theorem 2.1 relies on properties of the Gauss map.  Let $M$ be 
the image of a conformal minimal immersion $h:D \rightarrow \bfR^3$.  
The Gauss map $G: D \to S^2$ for the conformal 
minimal immersion $h:D \rightarrow M \subseteq \bfR^3$ maps
each point $p \in D$ to the unit normal of $M$ at $h(p)$ 
(considered as a point in the standard unit sphere $S^2$).  
$G$ is a holomorphic map from the complex 
coordinate $z = x_1 + i x_2 \in D$ to $S^2$ with 
the standard complex structure.  Therefore, if $G$ is not constant, it must map 
open sets in the interior of $D$ to
open sets in $S^2$.  We now prove Theorem 2.1.  

\begin{proof}
Let $h: D \rightarrow M \subset \bfR^3$ be any Douglas-Rad\'{o} solution 
for $\alpha$, as in Theorem 1.1.  Thus 
$h$ is a $C^\infty$ harmonic conformal 
minimal immersion on $D \setminus \{
h^{-1}(\ell_m \cap \ell_1), 
h^{-1}(\ell_1 \cap \ell_2), ... , 
h^{-1}(\ell_{m-1} \cap \ell_m) \}$, and 
$G$ is holomorphic on this same set.  
Also, $G$ is well-defined
and continuous at the vertices 
$h^{-1}(\ell_m \cap \ell_1), 
h^{-1}(\ell_1 \cap \ell_2), ... , 
h^{-1}(\ell_{m-1} \cap \ell_m)$; 
in fact, the unit normal at $h^{-1}(\ell_i \cap \ell_{i+1})$ (resp. 
$h^{-1}(\ell_m \cap \ell_1)$) must be perpendicular to the plane containing 
$\ell_i \cup \ell_{i+1}$ (resp. $\ell_m \cup \ell_1$) 
(\cite{DHKW}, Section 8.3).  

We now give the proof in five steps.  

\n {\bf Step 1:} {\em $G(h^{-1}(\ell_2 \cup ... \cup
\ell_{m-2} \cup \ell_m)) \subset S^2 \cap \{x_3 \geq 0\}$.}

Since $\cal P$ is a bijection from 
$\ell_2 \cup ... \cup \ell_{m-2}$ to 
$\partial \mbox{Conv}(P) \cap \partial P$, $\ell_2 \cup ... \cup 
\ell_{m-2}$ is contained in the boundary of the convex hull of $\alpha$.  
Since $h(D)$ is contained in the convex hull of $\alpha$ (see, for example, 
\cite{DHKW}, Section 6.1), 
the boundary point maximum principle \cite{Scn} implies that 
$G$ is never horizontal on $\ell_2 \cup ... \cup \ell_{m-2}$, except 
possibly at corner points $\ell_i \cap \ell_{i+1}, i = 
1,...,m-2$.  However, for $i = 1,...,m-2$, neither $\ell_i$ nor $\ell_{i+1}$ is 
vertical, and $\rho_i$ and $\rho_{i+1}$ are not parallel, and 
${\cal P}(\ell_i) = \rho_i$ and ${\cal P}(\ell_{i+1}) = 
\rho_{i+1}$, hence the normal 
vector at $\ell_i \cap \ell_{i+1}$ is not horizontal.  
Thus $G(h^{-1}(\ell_2 \cup ... \cup \ell_{m-2})) \cap \{x_3 = 0\} = \emptyset $.  
We choose the orientation of $M$ so that $G(h^{-1}(\ell_2 \cup ... \cup
\ell_{m-2})) \subset S^2 \cap \{x_3 > 0\}$.  Since $\ell_m$ is
vertical, $G(h^{-1}(\ell_m)) \subseteq S^2 \cap \{x_3 = 0\}$.  This shows Step 1.  

\n {\bf Step 2:} {\em There exists a horizontal vector $\vec{v} \in S^2$ such 
that $\vec{v} \not\in G(D)$.}  

By the conditions on $P$ and $\alpha$, there exists a
horizontal vector $\vec{v} = (v_1,v_2,0)$ so that any plane perpendicular to
$\vec{v}$ intersects $\alpha$ in at most two components.  We can
choose $\vec{v}$ so that for any plane $H$ perpendicular to $\vec{v}$ satisfying 
$H \cap \alpha \neq \emptyset$, one 
component of $H \cap \alpha$ is a single point, and the
other component is either empty or a
single point or $\ell_m$.  

We claim that there cannot be any point in the interior of
$M$ with normal $\pm \vec{v}$.  Suppose there is such a point 
$p \in \mbox{Int}(M)$.  Let $S(T_p(M)) = \{z \in D \; | \; h(z) \in
T_p(M)\}$.  By Lemma 2.2, $S(T_p(M))$ is a
plane embedded graph in $D$, and each vertex of $S(T_p(M))$ contained
in Int($D$) is connected to at least four edges.  
Note that $S(T_p(M)) \cap \mbox{Int}(D)$ 
has at least one vertex, at $h^{-1}(p)$.  

Since $T_p(M) \cap \alpha$ has at most two components, and the map 
$h|_{\partial D} : \partial D \rightarrow \alpha$ is bijective and 
continuous (Theorem~1.1), we know that $S(T_p(M)) \cap \partial D$ 
also has at most two components.  

It follows from elementary graph theory that 
$S(T_p(M))$ contains a closed loop $\beta$.  $h(\beta) \subseteq 
M \cap T_p(M)$ and $h(\beta)$ must bound a subdisk of $M$; and this subdisk 
must be contained in $T_p(M)$, since $h$ is 
a harmonic map.  Thus $M \subset T_p(M)$, since 
$h$ is harmonic.  But $\partial M = \alpha \not\subset T_p(M)$, 
thus there cannot be any point in Int($M$) with normal $\pm \vec{v}$.  

We now claim that there is at most one point in $\alpha$ 
with normal $\pm \vec{v}$.  Suppose there are two distinct points $p,q \in
\alpha$ with normal $\pm \vec{v}$.  By our choice of $\vec{v}$, 
the points $p,q$ must be contained in the interior of $\ell_m$.  
Let $H$ be the plane perpendicular to $\vec{v}$ and containing $\ell_m$.  
Let $S(H) = \{z \in D \; | \; h(z) \in H\}$.  
$S(H) \cap \partial D$ has two components, one of which is 
$h^{-1}(\ell_m)$.  There are (at 
least) two edges of $S(H)$ in Int($D$) with endpoints in 
$h^{-1}(\ell_m)$, meeting 
$h^{-1}(\ell_m)$ at the vertices $h^{-1}(p)$ and $h^{-1}(q)$.  
Thus again we see $S(H)$ must contain a closed loop, and 
we have a contradiction.

Therefore either $\vec{v} \not\in G(D)$ or $-\vec{v} \not\in G(D)$.  
Changing $\vec{v}$ to $-\vec{v}$ if necessary, we have 
$\vec{v} \not\in G(D)$.  This shows Step 2.  

\n {\bf Step 3:} {\em $G(D) \subset S^2 \cap \{x_3 \geq 0\}$, and for 
any point $z \in $Int($D$) there is an open neighborhood 
$U \subset $Int($D$) of $z$ so that 
$h(U)$ is a graph over $\{ x_3 = 0 \}$.}  

Since $h(D)$ is contained in the convex hull of $\alpha$, and since 
$\ell_{m-1}$ and $\ell_1$ are contained in the boundary of
the convex hull of $\alpha$, $G(h^{-1}(\ell_{m-1})) \subseteq \sigma_{m-1}$ and 
$G(h^{-1}(\ell_1)) \subseteq \sigma_1$, where $\sigma_{m-1}$ and 
$\sigma_1$ are 180 degree arcs of great circles in $S^2$ from $(0,0,1)$ to 
$(0,0,-1)$.  Furthermore,
the boundary point maximum principle \cite{Scn} implies 
$(0,0,-1) \not\in G(h^{-1}(\ell_{m-1} \cup \ell_{1}))$.  

We saw in Step 1 that $G(h^{-1}(\ell_2 \cup ... \cup
\ell_{m-2} \cup \ell_m)) \subset S^2 \cap \{x_3 \geq 0\}$ and 
$G(h^{-1}(\ell_m)) \subseteq S^2 \cap \{x_3 = 0\}$.  This and the preceding 
paragraph imply that 
\[ {\cal X} = (S^2 \cap \{x_3 \leq 0\}) \setminus \{ G(h^{-1}(\ell_{m-1} \cup 
\ell_{1})) \} \] 
is a connected set.  

We will show that 
\[ G(D) \subset S^2 \cap \{x_3 \geq 0\} \; . \] Suppose 
$G(D) \not\subset S^2 \cap \{x_3 \geq 0\}$.  Then, by Step 1, there is some 
point $p \in \mbox{Int}(D) \cup h^{-1}(\ell_1) \cup h^{-1}(\ell_{m-1})$ such 
that $G(p) \in S^2 \cap \{ x_3 < 0 \}$.  Hence some open neighborhood $U$ of 
$p$ in $D$ satisfies $G(U) \subset S^2 \cap \{ x_3 < 0 \}$.  Since the Gauss 
map $G$ maps open sets to open sets, we have 
\[ \mbox{Int}({\cal X}) \cap G(D) \neq \emptyset \; . \]  Since 
${\cal X} \cap G(D)$ is 
both open and closed in $\cal X$, and since $\cal X$ is connected, we have 
\[ {\cal X} \cap G(D) = {\cal X} \; . \]  
(${\cal X} \cap G(D)$ is closed in $\cal X$, since $G$ is holomorphic, and 
$D$ is closed; ${\cal X} \cap G(D)$ is open in $\cal X$, since $G$ is 
holomorphic and so $\partial G(D) \subset G(\partial D)$, and also since 
$G(\partial D) \cap \mbox{Int}({\cal X}) = \emptyset$.)  
Therefore ${\cal X} \subset G(D)$ and so $S^2 \cap \{ x_3 \leq 0 \} \subset 
G(D)$, contradicting Step 2.  
We conclude that $G(D) \subset S^2 \cap \{x_3 \geq 0\}$.  

Finally, by the holomorphicity of $G$, 
$G(\mbox{Int}(D)) \subset S^2 \cap \{x_3 > 0\}$.  
Thus at each point in Int($M$), $M$ is locally a graph over the $x_1x_2$-plane.  
This shows Step 3.  

\n {\bf Step 4:} {\em $M$ is embedded and ${\cal P}($Int$(M)) 
\subseteq $ Int($P$).}  

Let $\partial {\cal P}(M)$ be the boundary of ${\cal P}(M)$ in $\{x_3 = 0 \}$.  
Suppose there exists a point $p \in $ Int($M$) such that ${\cal P}(p) \not\in 
$ Int($P$), then $(\partial {\cal P}(M)) \setminus P$ is 
not empty.  Let $\ell$ be a vertical line intersecting 
$(\partial {\cal P}(M)) \setminus P$.  The line $\ell$ must make a
tangential intersection with some point $q \in $ Int($M$).  Thus 
$T_q(M)$ is a vertical tangent plane.  This
contradicts Step 3, hence ${\cal P}$(Int($M$)) $\subseteq$ Int($P$).  
We conclude that $M$ is embedded at its boundary.  
Thus, by Lemma 2.1, $M$ is embedded.  This shows Step 4.  

\n {\bf Step 5:} {\em  
$h|_{\mbox{Int}(D)}$ is a graph over Int($P$), and is the unique 
Douglas-Rad\'{o} solution with boundary $\alpha$.}  

The arguments in the next two paragraphs are similar to the proof of 
Theorem 1 in \cite{Scn}, except that our projection domain $P$ is 
not convex, and we use a family of translations instead of Schoen's family 
of reflections.  Hence we only outline the arguments here.  

First we show $h:$Int($D$)$\rightarrow $Int($M$) is a graph over 
Int($P$).  Let 
$h_\lambda:D \rightarrow M$ be defined by $h_\lambda(p) = h(p) +
(0,0,\lambda)$.  Choose $\lambda_0 \geq 0$ to be the smallest value so that 
for any $\lambda \geq \lambda_0$, $h_\lambda(\mbox{Int}(D))$ and 
$h(\mbox{Int}(D))$ have no points of transverse intersection.  
If $\lambda_0 > 0$, then $h(D)$ and $h_\lambda(D)$ must violate the maximum
principle \cite{Scn}, either at an interior point or at a boundary
point.  Thus $\lambda_0 = 0$, which implies that Int($M$) is a graph.
(Note that we are using $h_\lambda(\mbox{Int}(D))$ and
$h(\mbox{Int}(D))$ to define $\lambda_0$, and we are not using
$h_\lambda(D)$ and
$h(D)$.  This distinction is important, as an intersection of
$h_\lambda(D)$ and
$h(D)$ at a point in $h_\lambda(\partial D)\cap  h(\partial D)$ does
not necessarily constitute a contradiction to the maximum principle.)

Finally, suppose there exist two Douglas-Rad\'{o} solutions 
$h:D \rightarrow \bfR^3$ and $g:D \rightarrow \bfR^3$.  As we have
shown, they must both be embedded graphs over $P$.  
Let 
$g_\lambda(p) = g(p) + (0,0,\lambda)$.  Choose $\lambda_0 \geq 0$ to be the
smallest value so that
for any $\lambda \geq \lambda_0$, $g_\lambda(\mbox{Int}(D))$ and
$h(\mbox{Int}(D))$ have no points of transverse intersection.  
If $\lambda_0 > 0$, the maximum principle is violated.  Thus
$\lambda_0 = 0$, which implies that $g(D)$ lies above $h(D)$.
Similarly, $h(D)$ lies above $g(D)$.  Therefore $g(D) = h(D)$, and the
Douglas-Rad\'{o} solution is unique.  This shows Step 5.  
\end{proof}

\section{The Conjugate Surface Construction}

The Weierstrass representation is a principal tool used for the
construction of minimal surfaces in $\bfR^3$.  
Given a compact Riemann surface $\Sigma$, a set of points $\{p_j\}$ in
$\Sigma$, 
a meromorphic function $g:\Sigma\setminus\{p_j\} \goesto \bfC$,
and a holomorphic one-form $\omega$ on $\Sigma \setminus \{ p_j \}$,
the mapping $X:\Sigma \setminus \{ p_j \} \goesto \bfR^3$ defined by 
\begin{equation}
        X(z) = \Re \int_p^z 
                \left( \frac{1}{2}(g^{-1} - g)\omega,
                \frac{i}{2}(g^{-1} + g)\omega,
                \omega \right)
\end{equation}
is a conformal minimal immersion, where
$p \in \Sigma$ is fixed.  
$X$ is regular away from poles and zeroes of $g$ provided $\omega$ is
nonzero there, and $X$ is regular 
at a pole or zero of $g$ of order $m$ provided $\omega$ has a zero 
there of order $m$.  
For $X$ to be well-defined on $\Sigma \setminus \{p_j\}$, we must 
have 
\begin{equation}
        \Re \oint_{\gamma} 
                \left( \frac{1}{2}(g^{-1} - g)\omega,
                \frac{i}{2}(g^{-1} + g)\omega,
                \omega \right) = 0
\end{equation}
for any representative $\gamma$ of any
non-trivial homotopy class.  

 The Riemann surface $\Sigma \setminus \{p_j\}$, meromorphic function $g$, and
one-form $\omega$ are referred to as the
{\em Weierstrass data}.
 Here $g$ is the Gauss map $G$ of $X$ composed with
stereographic projection to the complex plane.  

The conjugate surface $X^\prime$ of $X$ is the minimal surface with
the
same underlying Riemann surface $\Sigma \setminus
\{p_j\}$, and the same meromorphic function $g$, but with 
holomorphic one-form $i\omega$.  (Note that 
$(X^\prime)^\prime = - X$.)  
The parametrization $X^\prime(p)$ may only be
well-defined on a covering of $\Sigma \setminus
\{p_j\}$, since equation (3.2) can hold for the Weierstrass data 
$\{g,\omega\}$ on $\Sigma \setminus
\{p_j\}$ without 
holding for the Weierstrass data $\{g,i\omega\}$ on $\Sigma \setminus 
\{p_j\}$.  

Thus we have the maps $z \rightarrow X(z)$ and $z \rightarrow
X^\prime(z)$ from simply connected domains of
$\Sigma \setminus \{p_j\}$ to $X$ and
$X^\prime$, respectively.  This induces a covering map 
$\phi: X^\prime(z) \rightarrow z \rightarrow X(z)$, 
the {\em conjugate map}, from $X^\prime$ to $X$.  
The following lemma is proven in \cite{Ka1}, \cite{Ka3}, \cite{Ka4}.  

\begin{lemma}
The conjugate map $\phi$ has the following properties:
\newcounter{num}
\begin{list}%
{\arabic{num})}{\usecounter{num}\setlength{\rightmargin}{\leftmargin}}

\item $\phi$ is an isometry;

\item $\phi$ preserves the Gauss map $G$;

\item $\phi$ maps planar principal curves in $X^\prime$ to planar
asymptotic 
curves in $X$, and maps planar asymptotic curves in $X^\prime$ to
planar principal curves in $X$; that is to say, $\phi$ maps
non-straight planar geodesics to straight lines, and vice versa.

\end{list}
\end{lemma}



It follows from the second and third properties of $\phi$ that a
planar
geodesic in $X^\prime$ contained in a plane $H$ 
is mapped by $\phi$ to a line in $X$ that is perpendicular to $H$.

\section{Limit Surface Lemma and Period Removal Lemma}

We use Lemma 4.1 to produce 
stable noncompact embedded minimal 
surfaces from compact embedded least-area disks.  
It is a slight variation of a lemma in \cite{Ro}, 
and the proof in \cite{Ro} applies to this case as well.  

\begin{lemma}
Let $\{\alpha_i\}_{i=1}^\infty$ be a sequence of compact Jordan contours 
in $\bfR^3$ so that the following conditions hold:
\begin{list}%
{\arabic{num})}{\usecounter{num}\setlength{\rightmargin}{\leftmargin}}
\item There is a positive integer $n$ so that,
for all $i$, $\alpha_i$ is a 
polygonal Jordan curve consisting of at most $n$ 
line segments;  
\item Each $\alpha_i$ bounds a least-area minimal disk $M_i$;  
\item $\{\alpha_i\}_{i=1}^\infty$ 
converges (in the topology of compact uniform 
convergence) to a noncompact polygonal curve $\alpha$ (not necessarily
connected), and 
$\alpha$ consists of a finite 
number of line segments, rays, and complete lines.  
\end{list}
Then a subsequence of $\{M_i\}_{i=1}^\infty$ converges to a nonempty 
stable minimal surface $M$ (possibly disconnected) with boundary 
$\alpha$.  Furthermore, if each $M_i$ is embedded, then $M$ is 
embedded.  
\end{lemma}

Lemma 4.2 is useful for solving a period problem at a catenoid end of a
minimal surface.  Consider a minimal surface $M$ (with boundary 
$\partial M$) with an end that is a 180 degree 
arc of a helicoid end.  Denote a neighborhood of this end 
by $E$.  
Suppose that outside a compact ball in $\bfR^3$ the boundary $\partial 
E$ is a pair of straight rays $r_1, r_2$.  These two rays are
necessarily parallel and pointing in opposite directions.  
The conjugate surface $E^\prime$ of $E$ is a 
surface with a 180 degree arc of a catenoid end 
that, outside a compact ball 
in $\bfR^3$, is bounded by two infinite planar geodesics 
$s_1, s_2$ asymptotic to catenaries.  The curves $s_1, s_2$ 
lie in parallel planes, and these planes are perpendicular to $r_1$
and $r_2$.  For this situation, we 
have the following lemma.  A proof can be found in \cite{Ro}.  

\begin{lemma}
The two planar geodesics $s_1, s_2 \subset \partial E^\prime$ lie 
in the same plane if and only if the plane containing the 
two conjugate straight boundary 
rays $r_1, r_2 \subset \partial E$ is parallel to the normal vector at the
helicoid end of $E$.  
\end{lemma}

\section{Complete Minimal Surfaces}

The conjugate surface construction 
method described in the introduction 
has been successful in many cases of 
minimal surfaces $M$ of the following type:

\begin{itemize}
\item $M$ has catenoid ends;

\item each end is invariant under some plane of reflective
symmetry of $M$;

\item the conjugate surface $M^\prime$ of the fundamental piece of $M$
is embedded;

\item all period problems that do not occur at an end of 
$M$ can be simultaneously removed by comparison arguments.  
\end{itemize}

Many previously known examples fit this description.  Among them are the
Jorge-Meeks $n$-oids \cite{JoMe}, the genus-1 $n$-oids \cite{BeRo}, 
the Platonoids \cite{Xu} \cite{Kat} \cite{UmYa}, the
higher-genus Platonoids \cite{BeRo}, the ${\cal A}{\cal W}_0(2n,w)$ surfaces
\cite{Ro}, the prismoids \cite{Ro}
\cite{Kat}, the higher-genus prismoids \cite{Ro}, and the
Jorge-Meeks fence \cite{Ro}.  

Some of the examples below have been shown to exist by other methods.  
Wohlgemuth \cite{Wo} has made similar periodic examples, by adding handles 
to a catenoid.  For his examples, he constructs the Weierstrass 
data.  In \cite{Ka3}, Karcher shows how 
to deal with examples similar to the first three examples below, by
directly constructing the Weierstrass data.  One 
might also be able to construct the last three examples using Weierstrass
data, using the methods of \cite{Ka1}, \cite{Kat}, and \cite{Wo}.

The purpose of the examples here is 
to demonstrate that we can apply the conjugate surface
construction to cases where it couldn't be applied before, by 
using Theorem 2.1.  

\begin{example}
{\em 
Choose any real number $w>0$ and any integer $n \geq 3$.  
For each positive integer $i$, let $\alpha_i$ be 
the polygonal curve with line segments from $(-1,w,0)$ to $(i,w,0)$ to 
$(i,0,i)$ to $(0,0,i)$ to $(0,0,-i)$ to 
$(-1,-i\tan(\frac{\pi}{n}),-i)$ and back to $(-1,w,0)$.
Let $\partial P_i$ be the 5-gon in the $x_1x_2$-plane such that 
${\cal P}(\alpha_i) = \partial P_i$.  
By Theorem 2.1, the Douglas-Rad\'{o} solution $M_i$ for 
$\alpha_i$ is unique and is a graph over $P_i$.  By Lemma
4.1, there exists an embedded limit surface $M$ for some subsequence of $i
\rightarrow \infty$.  (See Figure 2.)  

The boundary of $M$ consists of a ray with endpoint $(-1,w,0)$
pointing in the direction of the positive $x_1$-axis, a ray with
endpoint $(-1,w,0)$ containing the point
$(-1,w-\tan(\frac{\pi}{n}),-1)$, and the $x_3$-axis.  By construction,
${\cal P}(M) \subseteq \lim_{i \rightarrow \infty} P_i$.  Since each $M_i$ is a graph
over $P_i$, we may also conclude that the image of the Gauss map
$G$ on $M$, $G(M) \subseteq S^2 \cap \{x_3 \geq 0\}$.  (Here we are making a 
convenient 
abuse of notation by considering $G$ to be defined directly on the minimal 
surface, rather than on some immersion of the surface.)  
Since each $M_i$ is a disk, $M$ is either a single 
simply-connected surface, or the union $M = M_A \cup M_B$ of 
two disjoint simply-connected surfaces $M_A$ and $M_B$, with $\partial M_A$ being 
the two rays extending from $(-1,w,0)$, and $\partial M_B$ being 
the $x_3$-axis.  However the second case $M = M_A \cup M_B$ is not 
possible, and we defer to the Appendix for a proof of this.  Thus $M$ is 
connected.  The argument in the Appendix also shows that $M$ has 
finite total curvature.  

\begin{figure}
        \hspace{1.0in}
        \epsfxsize=1.1in
        \epsffile{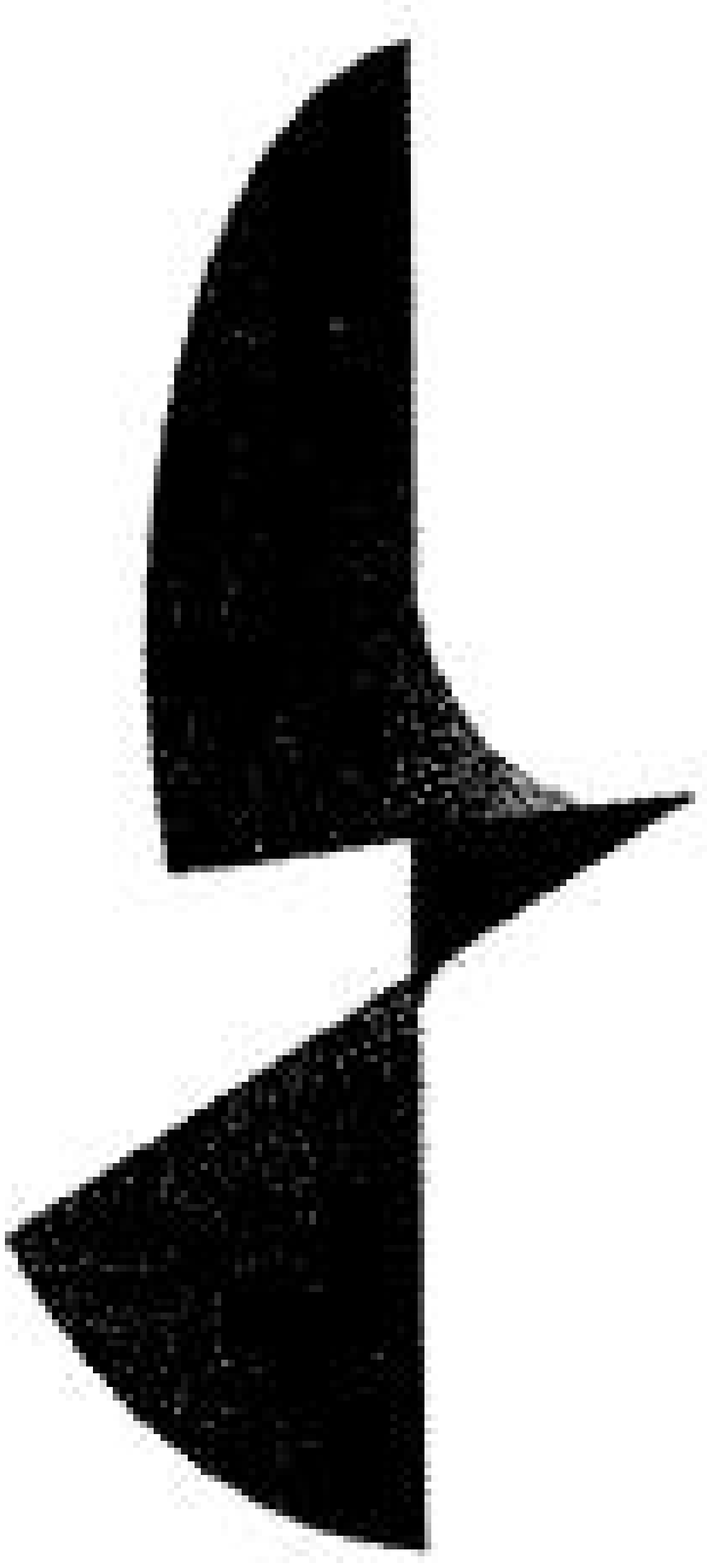}
        \hfill
        \epsfxsize=2.0in
        \epsffile{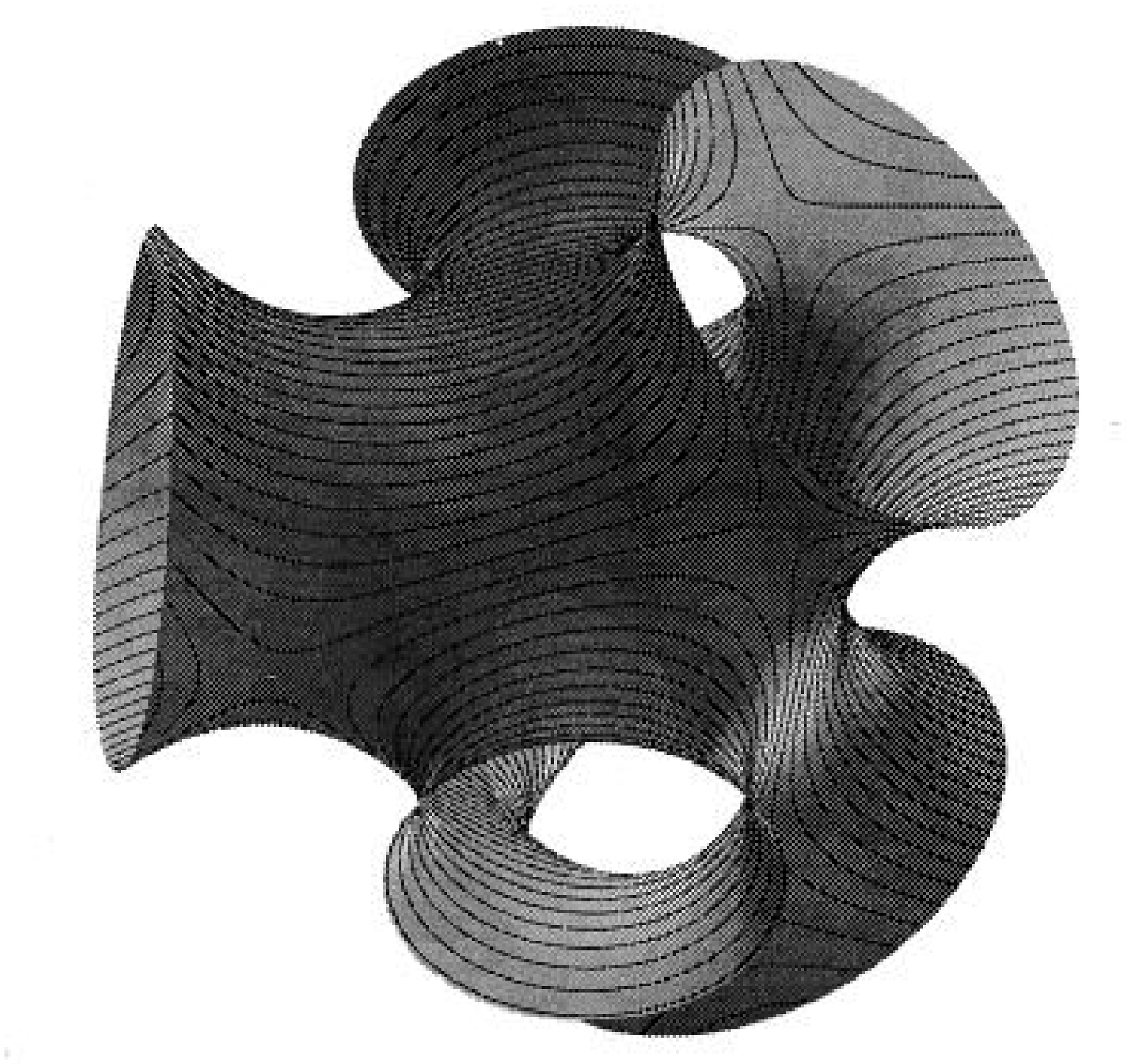}
        \hspace{1.0in}
        \caption{The limit surface $M$ described in Example~5.1 and
          the resulting complete minimal surface, with $n=3$.}
\end{figure}

The conjugate surface $M^\prime$ to $M$ is bounded
by three planar geodesics, none of which lie in parallel planes.
$M^\prime$ has one end which is a 90 degree arc of a catenoid end, and
another end which is a 180/$n$ degree arc of a catenoid end.  Thus 
we can extend $M^\prime$ using Theorem 1.2 to a complete
minimal surface.  This surface has no period problems, and therefore is
nonperiodic and of finite total curvature.  It consists of 4$n$ copies
of $M^\prime$, and has $n+2$ catenoid ends.  
Amongst the $n+2$ ends, $n$ of them have equal weight, and the other 
two have equal weight.  By a homothety of $\bfR^3$ if necessary, we may
assume that $n$ ends have weight 1 and two ends have weight $r=cw$ for
some positive constant $c$ depending only on $n$.
This surface is known to exist by other methods \cite{Kat},
\cite{Xu}.  It was also proven to exist by the conjugate surface
construction in Theorem 1.3 of \cite{Ro}.  However, in \cite{Ro} the
additional assumption was made that $r$ is larger than some given 
positive constant.  Here, due to Theorem 2.1, we can show
existence of the surface for any $r>0$.  (See Figure 2.)

We remark that this example can also be constructed, as above, for the 
case $n=2$.  For $n=2$, we need only replace the vertex 
$(-1,-i\tan(\frac{\pi}{n}),-i)$ of $\alpha_i$ with $(-1,-i,0)$ instead.  
}
\end{example}

\begin{example}
{\em 
Choose any integer $n \geq 2$, and any real numbers $w>0$ and 
$s > w/\sin(\frac{\pi}{n})$.  
For each positive integer $i$, let $\alpha_i$ be the polygonal curve 
from $(-1,w,0)$ to 
$(i,w,0)$ to $(i,0,i)$ to $(0,0,i)$ to $(0,0,-i)$ to 
$(i(s \cdot \sin(\frac{\pi}{n})-w),-i,-i)$ to 
$(\frac{s \cdot \sin(\frac{\pi}{n})-w}{i}-1,w-s \cdot \sin(\frac{\pi}{n})-
\frac{1}{i},i)$ to 
$(-1,w-s \cdot \sin(\frac{\pi}{n}),-s \cdot \cos(\frac{\pi}{n}))$ and back to 
$(-1,w,0)$.  
Let $\partial P_i$ be the 7-gon in the $x_1x_2$-plane so that ${\cal P}
(\alpha_i) = \partial P_i$.  
By Theorem 2.1, the Douglas-Rad\'{o} solution $M_i$ for 
$\alpha_i$ is unique and is a graph over $P_i$.  Again
some subsequence converges to an embedded surface $M$ as $i \rightarrow 
\infty$.  As in the last example, one can argue 
that $M$ is simply connected and of finite total curvature (see the Appendix).  


\begin{figure}
        \hfill
        \epsfxsize=3.1in
        \epsffile{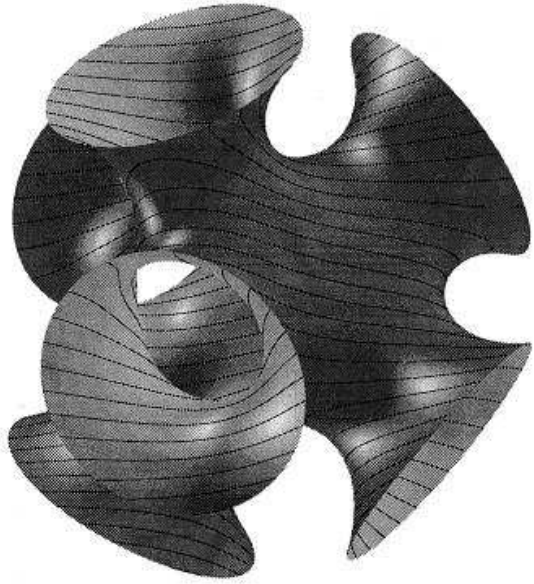}
        \hspace{1.4in}
        \caption{A prismoid with two layers of ends and with n=3.  
        The 3-layered prismoid is similar, but has an additional $n$
        ends along the horizontal plane of symmetry.}
\end{figure}

Let $M^\prime$ be the conjugate surface to $M$.  
There are two planar geodesics in $\partial M^\prime$ that lie in
parallel planes.  In order to extend using Theorem 1.2 to a complete
minimal surface with finite total curvature, these two parallel planes
must be the same plane.  Thus there is one
period problem to solve at a catenoid end of the surface.  
By Lemma 4.2, these two parallel planes are equal if the vertical ray
and complete vertical line in $\partial M$ lie in a common plane perpendicular
to the 180 degree helicoid end of $M$.  This is the case, since we
constructed $\alpha_i$ so that this would be so.  Thus the
period problem 
is solved, and using Theorem 1.2, $M^\prime$ extends 
to an immersed complete minimal surface of finite
total curvature with catenoid ends.  The complete surface consists of
4$n$ copies of $M^\prime$, and is a
``prismoid'' with 3 layers of ends, and its symmetry
group is $D_n \times \bfZ_2$.  It may be placed in $\bfR^3$ so that it 
has $n$ ends with horizontal normal
vectors, all of equal weight, and has $n$
ends with normal vectors pointing upward making an angle $\theta$ with
a horizontal plane, and has $n$
ends with normal vectors pointing downward making the same 
angle $\theta$ with a horizontal plane, for any $\theta \in
(0,\pi/2)$.  All
of these last 2$n$ ends have equal weight.  The ratio between the
weight of the first $n$ ends and the weight of the last 2$n$ ends 
can be any positive value.  Thus, for each $n$, 
we have a two-parameter family of 
these surfaces.  This example 
has been shown to exist by a different method in \cite{Kat}.  (See 
Figure 3.) 
}
\end{example}

And we can produce examples that were previously unknown, as
in the examples below.

\begin{example}
{\em 
Choose any real numbers $w>0$ and $\lambda>0$, and choose
any integer $n \geq 2$.  Choose $y > \cot(\frac{\pi}{n})$ to 
be the unique value so that the distance from the point $(1,y,0)$ to 
the plane 
$\{x_2 = \cot(\frac{\pi}{n}) \cdot x_1 \}$ is $w$.  For each positive 
integer $i$, let 
$\alpha_i$ be the polygonal curve from 
$(0,0,-\lambda)$ to 
$(-1/i,1/i,\lambda))$ to $(0,i^2,\lambda+i)$ to 
$(1,y,\lambda+i)$ to $(1,y,-\lambda-i)$ to 
$(i^2 \sin(\frac{\pi}{n}),i^2 \cos(\frac{\pi}{n}),-\lambda-i)$ and back 
to $(0,0,-\lambda)$.  
Let $\partial P_i$ be the 5-gon in the $x_1x_2$-plane so that 
${\cal P}(\alpha_i) = \partial P_i$.  
By Theorem 2.1, for all $i$ sufficiently large, 
the Douglas-Rad\'{o} solution for 
$\alpha_i$ is unique and is a graph over $P_i$.  
Again, some subesquence converges to a limit surface $M$.  
As in the previous examples, we can show that $M$ is connected (see the 
Appendix).  The conjugate 
$M^\prime$ of $M$ can be extended by reflection to a complete minimal
surface in $\bfR^3$.  In this case there is one period problem that is
not at an end.  But here we do not solve the period problem,
as we wish to produce a periodic surface.  The resulting
surface is a periodic 
fence of ${\cal A}{\cal W}_0(2n,w)$ surfaces.  The ${\cal A}{\cal W}_0(2n,w)$
surfaces are described in \cite{Ro}, \cite{Kat}, and they essentially
look like Jorge-Meeks surfaces with 2$n$ ends, but the ends have
weights that alternate between two positive values.  
We have a 1-parameter family of these ${\cal A}{\cal W}_0(2n,w)$ fences,
given by the parameter $\lambda > 0$.  (See Figure 4.) 

\begin{figure}
        \hspace{0.0in}
        \epsfxsize=1.5in
        \epsffile{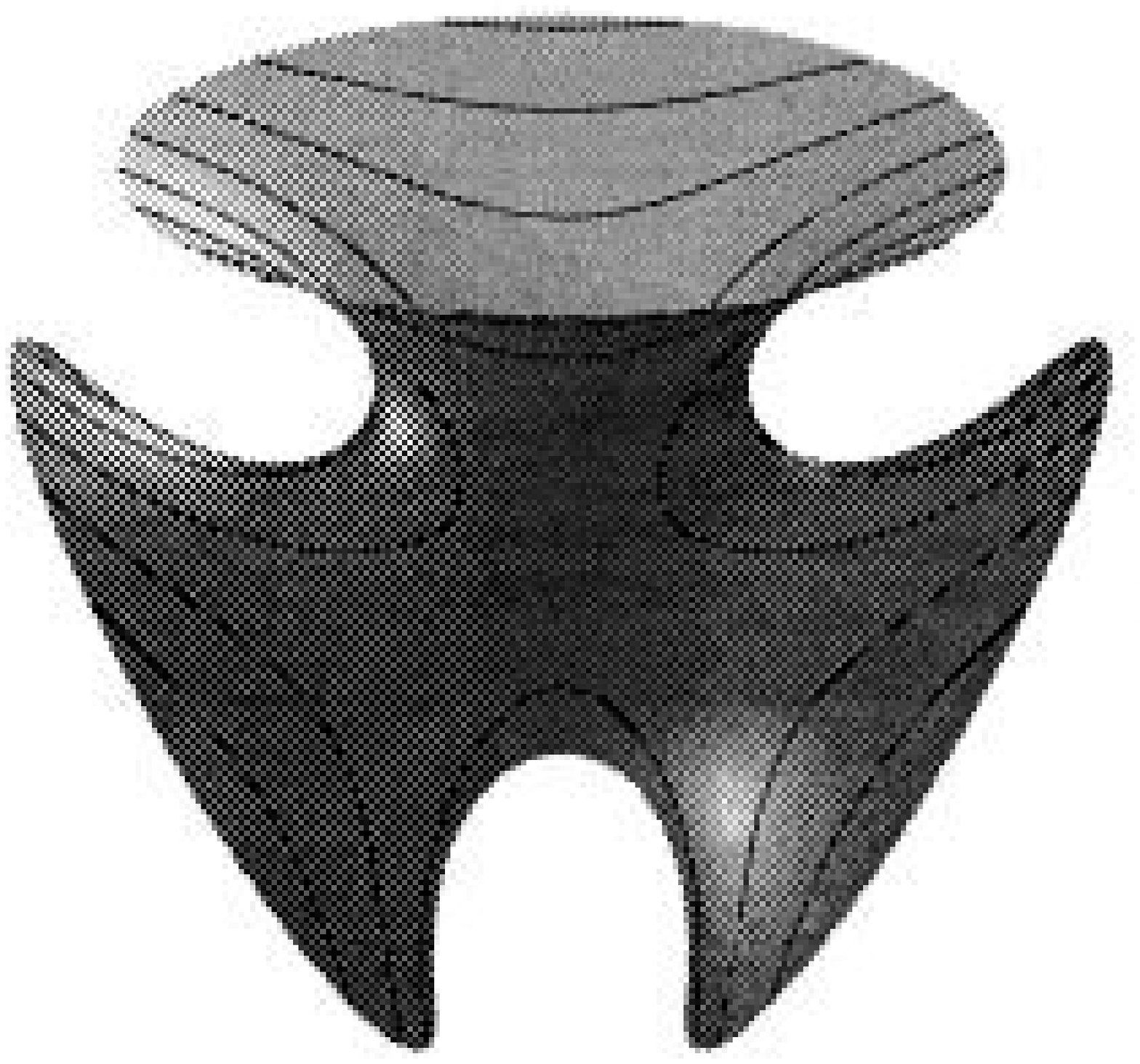}
        \hspace{0.6in}
        \epsfxsize=1.5in
        \epsffile{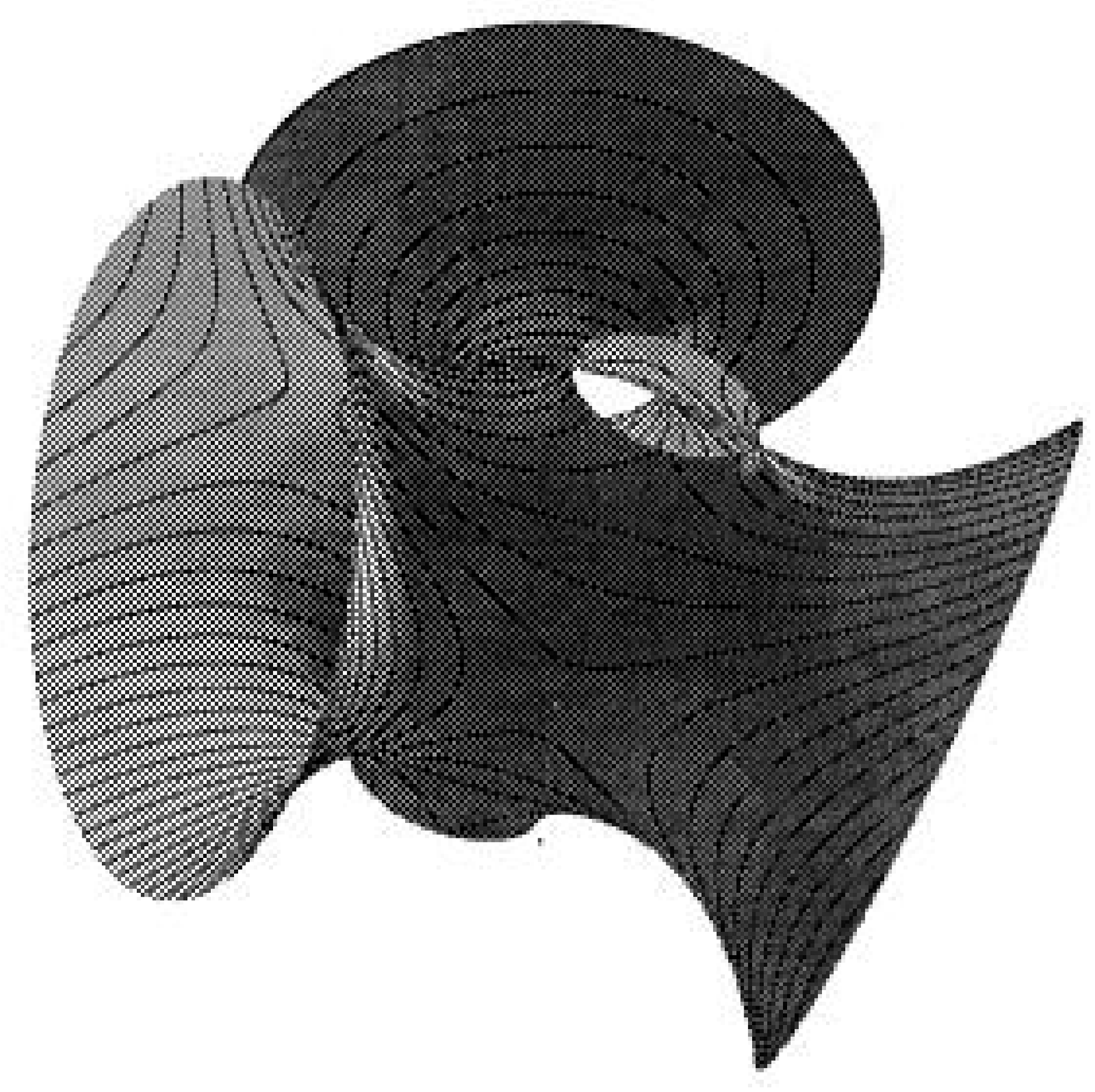}
        \hfill
        \epsfxsize=1.5in
        \epsffile{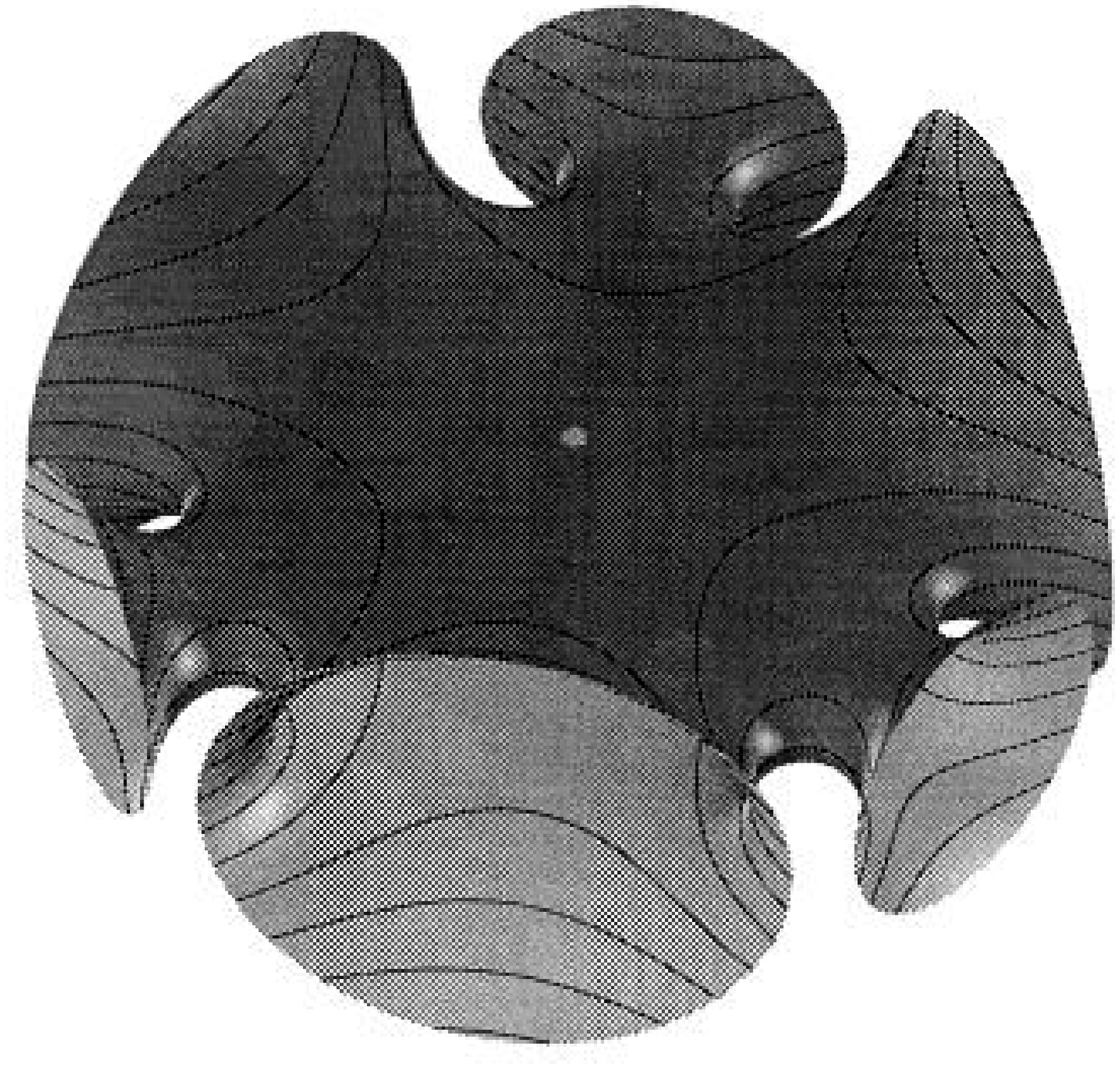}
        \hspace{0.3in}
        \caption{The Jorge-Meeks surface for $n=3$, the Jorge-Meeks
          fence for $n=3$, and a ${\cal A}{\cal W}_0(2n,w)$ surface
          for $n=3$.  
          Example 5.3 shows that ${\cal A}{\cal W}_0(2n,w)$ surfaces
          can be put together to make a periodic fence, just as
          Jorge-Meeks surfaces can be put together to make a
          Jorge-Meeks fence.} 
\end{figure}
}
\end{example}

\begin{example}
{\em
Let $n \geq 2$ be any integer, and let $\theta$ and $w$ be any real numbers 
such
that $0<\theta < \frac{\pi}{n}$ and 
$0<w<\frac{\sin(\frac{\pi}{n})}{\cos(\theta)}$.  For each positive
integer $i$, let $\alpha_i$ be the polygonal curve from 
$(0,0,0)$ to $((1-\frac{1}{i})\sin(\frac{\pi}{n}),
-(1-\frac{1}{i}-\frac{1}{i^2})\cos(\frac{\pi}{n}),0)$ to 
$(\sin(\frac{\pi}{n}),-\cos(\frac{\pi}{n}),-i)$ to 
$(i\sin(\theta), -i\cos(\theta), i)$ to 
$(\sin(\frac{\pi}{n})-w \cos(\theta), -\cos(\frac{\pi}{n})-w 
\sin(\theta), i)$ to 
$(\sin(\frac{\pi}{n})-w \cos(\theta), -\cos(\frac{\pi}{n})-w 
\sin(\theta), -i)$ to 
$(\sin(\frac{\pi}{n})-w \cos(\theta), -i^2, -i)$ and back to 
$(0,0,0)$.  Let $\partial P_i$ be the 6-gon so that ${\cal P}(\alpha_i)
= \partial P_i$.  
By Theorem 2.1, for any $i$ sufficiently large, 
the Douglas-Rad\'{o} solution for 
$\alpha_i$ is unique and is a graph over $P_i$.  As in the previous
examples, we can create a complete minimal surface with catenoid ends.
There is one period problem at an end which is solved by Lemma 4.2.

The resulting surface has a circle of ends, all of which are symmetric
across the same plane of reflective symmetry.  There are 3$n$ ends.
Up to a homothety, we may assume that $n$ of the ends have weight 1,
and that the other 2$n$ ends have weight $r$.  We may choose $r$ to be
any positive number.  
As one travel around this circle of ends, the weights of the ends
follow a pattern of 1,$r$,$r$,1,$r$,$r$,...,1,$r$,$r$.  
The angle
between any two adjacent ends with different weights is 
$\theta$.  
The angle between any two adjacent ends both of weight $r$ is 
$2(\frac{\pi}{n}
- \theta)$.  
Thus, for each $n$, 
we have a 2-parameter family of these surfaces, with parameters
$\theta$ and $r$.  
}
\end{example}

\begin{example}
{\em
One can also produce a genus-1
counterpart to the last example, just as the genus-1 $n$-oid is a
genus-1 counterpart to the genus-0 Jorge-Meeks $n$-oid.  (See Figure
5.)  The author has verified that one can construct finite contours
$\alpha_i$ so that Theorem~2.1 can be applied to the genus-1 case as
well.  As before, we have a connected limit surface $M$ and a
conjugate fundamental piece $M^\prime$.  In this case $M^\prime$ has
two period problems.  
One of them is at a catenoid end and is solved by
Lemma 4.2.  The other is not at an end, and we can solve this by a
comparison argument using a portion of a helicoid.  We do not include
the comparison argument here, as it is similar to arguments in 
\cite{Ka4}, \cite{BeRo}, and \cite{Ro}.


\begin{figure}
        \hfill
        \epsfxsize=3.1in
        \epsffile{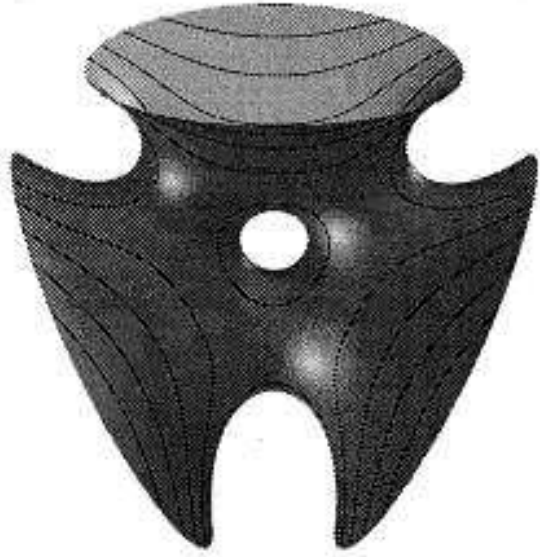}
        \hspace{1.4in}
        \caption{The genus-1 $n$-oid for $n=3$.}
        \label{g1tri}
\end{figure}

}
\end{example}

\section{Appendix}

In this Appendix, we will show that the limit surface $M$ is 
connected and of finite total 
curvature in each of Examples 5.1, 5.2, 5.3, 5.4, and 5.5.

{\bf For Example 5.1:}  In this example, the surfaces $M_i$ have boundaries 
$\partial M_i$ consisting of six line segments and six vertices.  
At five of the vertices the exterior angle is $\frac{\pi}{2}$ radians, 
and at the other vertex the exterior angle approaches
$\frac{\pi}{2} + \frac{\pi}{n}$ as $i \to \infty$.  The Gauss-Bonnet 
theorem then implies that 
\[ \int_{M_i} |K| dA = - \int_{M_i} K dA \to \pi +\frac{\pi}{n} 
\] as $i \to \infty$.  (Note that $dA$ is the area form on $M_i$ induced 
as a submanifold of $\bfR^3$, and that $|K| = -K$ on a minimal surface.)  
Thus the limit surface $M$ has finite total curvature at most 
$\pi + \frac{\pi}{n}$.  

Suppose the second case $M = M_A \cup M_B$ described in Example 5.1 
occurs; that is, suppose $M$ is the union of two simply-connected 
minimal surfaces $M_A$ and $M_B$.  $M_A$ is embedded, of finite 
total curvature, and is bounded by two rays as described in Example 5.1.  
Furthermore, ${\cal P}(M_A) \subseteq \lim_{i \to 
\infty} P_i$, since ${\cal P}(M) \subseteq \lim_{i \to 
\infty} P_i$.  We will show that such an $M_A$ cannot exist, 
deriving a contradiction that implies $M$ is connected.  

The conjugate surface $M_A^\prime$ of $M_A$ has 
boundary $\partial M_A^\prime$ consisting of two planar geodesics: 
one contained in the plane $\{x_1 = c_1 \}$ for some constant 
$c_1 \in \bfR$, the other contained in $\{x_3 + \tan(\frac{\pi}{n}) 
x_2 = c_2\}$ for some constant $c_2 \in \bfR$.  Let Ref$_1:\bfR^3 \to 
\bfR^3$ be reflection across 
the plane $\{x_1 = c_1 \}$, and let Ref$_2:\bfR^3 \to \bfR^3$ be reflection 
across the plane 
$\{x_3 + \tan(\frac{\pi}{n}) x_2 = c_2\}$.  Then 
$\hat{M}_A^\prime := M_A^\prime \cup \mbox{Ref}_1(M_A^\prime) \cup 
\mbox{Ref}_2(M_A^\prime \cup \mbox{Ref}_1(M_A^\prime))$ is a complete 
minimal surface with finite total curvature, and is 
simply connected with a single end.  By Theorem 9.5 of \cite{Os2}, 
the Gauss map $G$ extends continuously across the end of $\hat{M}_A^\prime$.  
Since the Gauss map is preserved by conjugation, 
$G$ extends continuously across the end of $M_A$ as well.  Thus the normal 
vector at the end of $M_A$ is well defined.  This normal vector 
must be perpendicular to both of the rays in $\partial M_A$, hence it is 
$(0,-1,\tan(\frac{\pi}{n}))$.  However, with this limiting normal vector 
at the end, it is clear that ${\cal P}(M_A) \not\subseteq \lim_{i \to 
\infty} P_i = \{(0,0,0)\} \cup \{(x_1,x_2,0) \, | \, x_1 \geq 0, x_2 \in 
(0,w] \} \cup \{(x_1,x_2,0) \, | \, x_1 \in [-1,0), x_2 \leq w \}$, a 
contradiction.  

{\bf For Example 5.2:}  Suppose the $M$ in Example 5.2 is not connected.  
It then consists of 
two disjoint embedded simply-connected minimal surfaces $M_A$ and 
$M_B$.  Let $M_A$ be the component bounded by the ray pointing in the 
direction of the positive $x_1$-axis with endpoint $(-1,w,0)$, the ray 
pointing in the direction of the positive $x_3$-axis with endpoint 
$(-1,w-s \sin(\frac{\pi}{n}),-s \cos(\frac{\pi}{n}))$, and the line 
segment from $(-1,w,0)$ to 
$(-1,w-s \sin(\frac{\pi}{n}),-s \cos(\frac{\pi}{n}))$.  Since 
${\cal P}(M_i) \subseteq P_i$, we have 
${\cal P}(M_A) \subseteq \lim_{i \to \infty} P_i$.  
Using the Gauss-Bonnet theorem just like for Example 5.1, we 
conclude that $M$ and $M_A$ have finite total curvature.  

The conjugate surface $M_A^\prime$ of $M_A$ is bounded by 
three planar geodesics, one of finite length, another of infinite 
length contained in the plane $\{x_1 = c_1\}$ for some constant 
$c_1 \in \bfR$, and 
the third of infinite length contained in the plane $\{x_3 = c_2\}$ for 
some constant $c_2 \in \bfR$.  Let Ref$_1:\bfR^3 \to \bfR^3$ be reflection 
across 
the plane $\{x_1 = c_1\}$, and let Ref$_2:\bfR^3 \to \bfR^3$ be reflection 
across the plane $\{x_3 = c_2\}$.  Then 
$\hat{M}_A^\prime := M_A^\prime \cup \mbox{Ref}_1(M_A^\prime) \cup 
\mbox{Ref}_2(M_A^\prime \cup \mbox{Ref}_1(M_A^\prime))$ is an 
annular minimal surface of finite total curvature with a single 
compact boundary loop and a single end.  

Unlike the case of Example 5.1, in this case the boundary $\partial 
\hat{M}_A^\prime \neq \emptyset$; however, $\partial \hat{M}_A^\prime$ is a 
compact loop and hence we may still apply 
Theorem 9.5 of \cite{Os2} to conclude that 
the Gauss map $G$ extends across the end of $\hat{M}_A^\prime$.  Hence 
$G$ extends to the end of $M_A$.  As in the case of Example 5.1, 
we see that ${\cal P}(M_A) \not\subseteq \lim_{i \to 
\infty} P_i$.  This contradiction implies $M$ is connected.  

{\bf For Examples 5.3, 5.4, and 5.5:}  In 
these final three examples it can be shown, in the same way as 
for Example 5.2, that $M$ is connected and has finite total curvature.

\end{document}